\documentclass[12pt]{article}
\usepackage{stmaryrd}
\usepackage{amsfonts}
\usepackage{mathrsfs}
\usepackage{amsmath,amscd}
\usepackage{amssymb,amsmath}
\def\Z{\mathbb{Z}}
\def\C{\mathbb{C}}
\def\QED{\hfill$\Box$}

\numberwithin{equation}{section}
\newtheorem{theo}{Theorem}[section]
\newtheorem{defi}[theo]{Definition}
\newtheorem{coro}[theo]{Corollary}
\newtheorem{lemm}[theo]{Lemma}
\newtheorem{prop}[theo]{Proposition}

\newtheorem{rema}[theo]{Remark}

\begin{document}

\begin{center}
{\Large {\bf Some finite properties for  vertex operator superalgebras}}\\
\vspace{0.5cm}

Chongying Dong$^1$\footnote{Supported by NSF grants and a faculty research fund from the University of California at Santa Cruz.}
\\
Department of Mathematics, University of California, Santa Cruz
CA 95064\\
Jianzhi Han\\
Wu Wen-Tsun Key Laboratory of Mathematics, USTC\\ Chinese Academy
of Sciences, Hefei 230026, China
\end{center}
\begin{abstract} Vertex operator superalgebras are studied and various results on rational
Vertex operator superalgebras are obtained. In particular,
the vertex operator super subalgebras generated by the weight $\frac12$ and weight $1$ subspaces are determined.
It is also established that if the even part $V_{\bar 0}$ of a vertex operator superalgebra $V$ is rational, so is $V.$
\end{abstract}

\vskip10pt \noindent{\bf 1. \
Introduction}\setcounter{section}{1}\setcounter{theo}{0}\setcounter{equation}{0}
\vskip5pt

The vertex operator superalgebras which are  natural
generalizations of vertex operator algebras have been studied extensively in \cite{DZ2},
\cite{DZ1}, \cite{KW}, \cite{L2}, \cite{L3} and \cite{X}. In this
paper, we study certain finite properties of vertex operator superalgebras
following \cite{DLM4}, \cite{DM}, \cite{DM1} and \cite{Ma}.

A vertex operator superalgebra $V=V_{\bar 0}\oplus V_{\bar 1}$ has
even part $V_{\bar 0}$ and odd part $V_{\bar 1}$ where $V_{\bar
0}$ consists of vectors of integral weights and $V_{\bar 1}$
consists of vectors whose weights are half integers but not
integers. So there is a canonical automorphism $\sigma$ of $V$
acting on $V_{\bar i}$ as $(-1)^i$ and $V_{\bar 0}$ is a vertex
operator algebra which is also a fixed point subalgebra of $V.$ So
a better understanding of the relationship between representation
theories of $V$ and $V_{\bar 0}$ is definitely useful for the
study of orbifold theory (see \cite{DVVV} and \cite{DLM3}). Even
the orbiford theory for vertex operator algebra with order $2$
automorphism has not been understood fully.

Rationality which is an analogue of semisimplicity of associative
and Lie algebras is probably the most important concept in the
representation theory of vertex operator superalgebra. We first
establish that if $V_{\bar 0}$ is rational then $V$ is rational
although we believe that the rationalities of $V$ and $V_{\bar 0}$
are equivalent from the orbifold theory.  The main tool is the
associative algebras $A_{g,n}(V)$ for $n\in \frac12\Z_+$ which are
generalizations of $A_{g,n}(V)$ as introduced and studied in
\cite{DLM4} (also see \cite{Z}, \cite{KW}, \cite{DLM2} and
\cite{DLM0}) where $g$ is an automorphism of $V$ of finite order.
It is established that $V$ is $g$-rational if and only if
$A_{g,n}$ is a finite dimensional semisimple associative algebra
for large $n.$ This is the key result to prove the rationality of
$V$ from the rationality of $V_{\bar 0}.$ Another characterization
of rationality is given through Ext functor.

Our investigation next centers around the vertex operator super
subalgebras of $V$ generated by homogenous subspaces of small
weights. The vertex operator subalgebra generated by $V_{\frac12}$
is a holomorphic vertex operator superalgebra $U$ associated to an
infinite dimensional Clifford algebra built from a finite
dimensional vector space with a nondegenerate symmetric bilinear
form. This enables us to decompose $V$ as a tensor product
$U\otimes U^c$ where $U^c$ whose weight $\frac12$ subspace is zero
is the commutant of $U$ in $V$ (\cite{FZ},
\cite{LL}). Moreover, the module categories of $V$ and $U^c$ are
equivalent. To study $V_1$ we need to understand the algebraic
structure of $V_1$ first. Under the assumptions that $V$ is
rational or $\sigma$-rational together with $C_2$-cofiniteness we
are able to show that $V_1$ is a reductive Lie algebra using the
modular invariance results from \cite{DZ2} and \cite{Z} and the
fact that $E_2(\tau)$ is not modular. Also the rank of $V_1$ and
the dimension of $V_{\frac12}$ are controlled by the effective
central charge. Furthermore, for any simple Lie subalgebra
${\mathfrak g}$ of $V_1,$ the vertex operator subalgebra generated
by ${\mathfrak g}$ is isomorphic to the vertex operator algebra
$L(k,0)$ which is the integrable highest weight module for the
affine Kac-Moody algebra $\widehat{\mathfrak g}.$  We also give a
rational vertex operator subalgebra which is a tensor product of
affine vertex operator algebras and a lattice vertex operator
algebra and whose weight one subspace is exactly the $V_1.$

We should point out that most of results in this paper have been obtained in the case $V$ is a vertex operator algebra in \cite{DLM4}, \cite{DM}, \cite{DM1},  \cite{Ma}. So this paper can be regarded as a ``super'' analogues of results presented in \cite{DLM4}, \cite{DM}, \cite{DM1}, \cite{Ma}. The main ideas and the broad outlines also follow from these papers. A lot of arguments are omitted if they are the same
as in the case of vertex operator algebras. On the other hand, there is a new phenomenon in the super case.
Namely, either rationality together with $C_2$-cofiniteness or  $\sigma$-rationality together with $C_2$-cofiniteness implies that $V_1$ is reducitve. This gives a strong evidence that rationality, $\sigma$-rationality of $V$ and rationality of $V_{\bar 0}$ are equivalent. But we have no idea how to establish this.

This paper is organized as follows. We recall various notions of
twisted modules for a vertex operator supealgebra and
$g$-rationality for any automorphism of finite order from
\cite{FLM2}, \cite{Z}, \cite{DLM2} and \cite{DZ1} in Section 2. In
Section 3, we define a series of associative algebras $A_{g,n}(V)$
for a vertex operator superalgebra $V$ and $n\in \Z_+.$ We exhibit
how to use $A_{n}(V)$ to prove rationality of $V$ from the
rationality of $V_{\bar 0}.$ It is also shown that if $V$ is
$C_2$-cofinite or rational then $V$ is finitely generated and the
automorphism group $Aut(V)$ is an algebraic group. The Section 4
is devoted to the study of vertex operator super subalgebra
generated by $V_{\frac12}.$ In Section 5 we show that if $V$ is
rational or $\sigma$-rational together with $C_2$-cofiniteness
then the weight one subsapce $V_1$ is a reductive Lie algebra
whose rank is bounded above by the effective central charge
$\tilde c$. Consequently, $\dim V_{\frac12}$ is bounded above by
$2\tilde c+1.$ Section 6 deals with the vertex operator subalgebra
of $V$ generated by $V_1$ and related.

We make the assumption that the reader is
familiar with the theory of vertex operator algebras as presented in  \cite{B}, \cite{DL}, \cite{FLM2} and
\cite{LL}.

\vskip10pt \noindent {\bf 2. \
Basics}\setcounter{section}{2}\setcounter{equation}{0}

\vskip5pt  In this section we give the definition of vertex operator superalgebra  and several notions of modules (cf. \cite{DLM1}, \cite{DZ1} \cite{FFR}, \cite{FLM2}, \cite{L2}, \cite{Z}).

We first recall that a super vector space is $\Z_2$-graded vector space
$V=V_{\bar 0}\oplus V_{\bar 1}$. The elements in $V_{\bar 0}$
(resp., $V_{\bar 1}$) are called even (resp., odd). Let $\tilde v$
be 0 if $v\in V_{\bar 0}$, and 1 if $v\in V_{\bar 1}$.

\begin{defi} A vertex operator superalgebra (VOSA) is a $\frac 1
2\Z$-graded vector space $$V=\bigoplus_{n\in\frac 1 2\Z}V_n=V_{\bar
0}\oplus V_{\bar 1}$$ with $V_{\bar 0}=\sum_{n\in\Z}V_n$ and
$V_{\bar 1}=\sum_{n\in\Z}V_{n+\frac 1 2}$ satisfying all the axioms
in the definition of vertex operator algebra except that the Jocobi
identity is replaced by:
\begin{eqnarray*}
z_0^{-1}\!\!\!\!\!\!&\delta
&\!\!\!\!\!\big(\frac{z_1-z_2}{z_0}\big)Y(u,z_1)Y(v_2,z)-(-1)^{\tilde
u\tilde v}z_0^{-1}\delta\big(\frac{-z_2+z_1}{z_0}\big)Y(v,z_2)Y(u,z_1)\nonumber\\
\!\!&=\!\!\!\!\!\!&z_2^{-1}\delta\big(\frac{z_1-z_0}{z_2}\big)Y(Y(u,z_0)v,z_2).
\end{eqnarray*}
\end{defi}

Throughout the paper we always assume that $V$ is a vertex
operator superalgebra unless otherwise stated.

\begin{defi}\label{auto} An automorphism $g$ of a VOSA $V$ is a
linear automorphism of $V$ preserving the vacuum vector $\mathbf{1}$
and the conformal vector $\omega$ such that the actions of $g$ and
$Y(v,z)$ on $V$ are compatible in the sense that
$$gY(v,z)g^{-1}=Y(gv,z)$$ for $v\in V$.
\end{defi}

Denote by $Aut(V)$ the set consisting of all automorphisms of $V$.
Observe that any automorphism of $V$ commutes with $L(0)$ and
hence preserves each homogeneous subspace $V_n.$ As a consequence,
any automorphism preserves both $V_{\bar 0}$ and $V_{\bar 1}.$
There is a canonical automorphism $\sigma$ of $V$ with
$\sigma|{V_{\bar i}}=(-1)^i$ associated to the $\Z_2$-grading of
$V$.

Let $g\in Aut(V)$ with finite order $T$, then we can decompose $V$
into eigenspaces of $g:$
$$V=\bigoplus_{r=0}^{T-1} V^r$$
where $V^r=\{v\in V|gv=e^{\frac{-2\pi ir}{T}}v\}$.
\begin{defi} A weak $g$-twisted $V$-module $M$ is a $\Z_2$-graded vector space equipped with a
linear map
\begin{eqnarray*}
 & &V\rightarrow (\mathrm{End} M)[[z, z^{-1}]],\\
& &v\mapsto Y_M(v, z)=\sum_{n\in\frac 1 T\Z}v_nz_{-n-1},
\end{eqnarray*}
such that for all $u\in V^r(0\leq r\leq T-1)$, $v\in V$ and $w\in W$
the following conditions hold:
\begin{eqnarray*}
&Y_M(u,z)=\sum_{n\in\frac{r}{T}+\Z}u_nz^{-n-1};\\
&u_nw=0,\ \ for\  n>>0;\\
&Y_M(1, z)=\mathrm{Id}_M;
\end{eqnarray*}
\end{defi}

\begin{eqnarray*}\label{jocabim}
z_0^{-1}\!\!\!\!\!\!&\delta
&\!\!\!\!\!\big(\frac{z_1-z_2}{z_0}\big)Y_M(u,z_1)Y_M(v,z_2)-(-1)^{\tilde
u\tilde v}z_0^{-1}\delta\big(\frac{-z_2+z_1}{z_0}\big)Y_M(v,z_2)Y_M(u,z_1)\nonumber\\
\!\!&=\!\!\!\!\!\!&z_2^{-1}\delta\big(\frac{z_1-z_0}{z_2}\big)^{-r/T}\delta\big(\frac{z_1-z_0}{z_2}\big)Y_M(Y(u,z_0)v,z_2).
\end{eqnarray*}

\begin{defi} An  admissible $g$-twisted $V$-module
is a  weak $g$-twisted $V$-module $M$ which carries a
$\frac{1}{T}{\Z}_{+}$-grading
$$M=\oplus_{n\in\frac{1}{T}\Z_+}M(n)$$
satisfying
$$
v_{m}M(n)\subseteq M(n+\mathrm{wt} v-m-1)$$
for homogeneous $v\in V.$ \end{defi}

\begin{defi}\label{d2.4} An ordinary $g$-twisted  $V$-module is
a weak $g$-twisted  $V$-module
$$M=\bigoplus_{\lambda \in{\C}}M_{\lambda}$$
such that $\dim M_{\lambda}$ is finite and for fixed $\lambda,$ $M_{\frac n
T+\lambda}=0$ for all small enough integers $n$ where
 $M_{\lambda}=\{w\in M|L(0)w=\lambda w\}$ and $Y_M(\omega,z)=\sum_{n\in\Z}L(n)z^{-n-2}.$
\end{defi}

We say $V$ is {\em g-rational} if every admissible $g$-twisted
$V$-module is completely reducible, i.e. a direct sum of simple
admissible $g$-twisted $V$-modules. $V$ is {\em g-regular} if the
category of weak $g$-twisted $V$-modules is semisimple, namely,
every weak $g$-twisted $V$-module is a direct sum of irreducible
weak $g$-twisted $V$-modules. If $g=1,$ we have the definitions of
rationality and regularity for vertex operator superalgebras.

The following definitions are given for vertex operator algebras
in \cite{DM1} and \cite{Z} and we extend  these to vertex operator
superalgebras here.

 A vertex operator superalgebra $V$ is said to be of {\em CFT type} in case
that the $L(0)$-grading on $V$ has no negative weights, and if the
degree-zero homogeneous subspace $V_0$ is 1-dimensional, i.e.
$V=\bigoplus_{n\in\frac 1 2\Z_+}V_n$ and $V_0=\C \mathbf{1}$. We say
$V$ is  of {\em strong} CFT type if $V$ statisfies the further
condition $L(1)V_1=0$.  $V$ is said to be $C_2$-{\em cofinite} in
case that $C_2(V)$ has finite codimension in $V$, where $C_2(V)$ is
the subspace of $V$ linearly spanned by all elements of the form
$u_{-2}v$ for $u,v\in V$.

For convenience, let us introduce the term {\em strongly
g-rational}  for a simple vertex operator superalgebra $V$ which
satisfies the following conditions:
\begin{itemize}\parskip-3pt
\item[(1)] V is of strong CFT type;
\item[(2)] V is $C_2$-cofinite;
\item[(3)] V is  $g$-rational.
\end{itemize}

\begin{defi}\label{d25} A bilinear form $(\cdot,\cdot)$ on $V$-module $M$ is
said to be invariant \cite{FHL} if it satisfies the
following condition
$$(Y(a,z)u, v)=(u,
Y(e^{zL(1)}(-z^{-2})^{L(0)}a, z^{-1})v)\ \ \ \ for\ a\in V,\ u,v\in
M.$$
\end{defi}

It is proved in  \cite{L1} and \cite{X} that there exists a linear
isomorphism from the space of invariant bilinear forms on $V$ to
$\mathrm{Hom}_{\C}(V_0/L(1)V_1,\C)$. This implies that there is a
unique, up to multiplication by a nonzero scalar, nondegenerate
symmetric invariant bilinear form  on $V$ if $V$ is simple and of
strong CFT type.

\vskip10pt
 \noindent{\bf3.
Rationality}\setcounter{section}{3}\setcounter{theo}{0}\setcounter{equation}{0}
\vskip5pt

In this section we give a characterization of rationality of a
vertex operator superalgebra $V$ in terms of rationality of vertex
operator subalgebra $V_{\bar 0}$.  We will show that if  $V_{\bar
0}$ is rational then $V$ is rational. We certainly believe that
the converse is also true. That is,  if $V$ is rational  then
$V_{\bar 0}$ is also rational. This is similar to a well-known
conjecture in orbifold theory: Let $V$ be a rational vertex
operator algebra and  $g$ is an order 2 automorphism of $V.$ Then
the fixed point vertex operator subalgebra is also rational. We
will establish some other results on rationality. We also discuss
the generators of $V.$

The tool we use to prove the main result is the associative algebras $A_n(V)$ which is defined in \cite{DLM0} for vertex operator algebra. Let $V$ be vertex operator superalgebra. Let $O_n(V)$ be the subspace of $V$ linearly spanned by all
$L(-1)u+L(0)u$  and $u\circ_nv$ where for homogeneous $u\in V$ and
$v\in V$
$$u\circ_nv=\left\{ \begin{aligned}
          &{\rm Res\it}_{z}\frac{(1+z)^{{\rm
wt\it}u+n}}{z^{2n+2}}Y(u,z)v,
\quad \mathrm{if} \ u\in V_{\bar{0}}, \\
                  &{\rm Res\it}_z \frac{(1+z)^{{\rm
wt\it}u+n-\frac 1 2}}{z^{2n+1}}Y(u,z)v, \quad \mathrm{if} \ u\in
V_{\bar{1}}.
                          \end{aligned} \right.$$
Define another operation $\ast_n$ on $V$ by
$$u\ast_nv=\left\{ \begin{aligned}
          &\sum\limits_{m=0}^{n}(-1)^{m}{m+n\choose
n}\mathrm {Res}_{z}\frac{(1+z)^{\mathrm{wt}u+n}}{z^{n+m+1}}Y(u,z)v,
\quad
\mbox{if} \ u,v\in V_{\bar{0}}, \\
                  & 0,  \ \ \ \ \ \ \ \ \ \ \ \ \ \ \ \ \ \ \ \ \ \ \ \
                  \ \ \ \ \ \ \ \ \ \ \ \ \ \ \ \ \ \ \ \ \ \ \ \ \ \quad
\mbox{if} \ u\in V_{\bar{1}}\ \mbox{or}\ v\in V_{\bar{1}}.
                          \end{aligned} \right.$$
Set $A_n(V)=V/O_n(V).$ Then $A_0(V)$ is the $A(V)$ studied in \cite{KW}.
Let $M$ be a weak $V$-module. Define the
``$n$-th lowest weight vector" subspace of $M$ to be
$$\Omega_n(M)=\{w\in M\mid u_{\mbox{wt}u+n+i}w=0, u\in V, i\geq 0\}.$$

As in \cite{DLM0} we have the following results:
\begin{theo}\label{theo3.2}
(1) Suppose that $M$ is a weak $V$-module. Then $\Omega_n(M)$ is
an $A_n(V)$-module such that  $a$ acts as $o(a)$ for  $a\in
V_{\bar 0}$, where $o(a)$ is defined to be $a_{\mathrm {wt}a-1}$
for homogeneous $a\in V_{\bar 0}$ and extends it linearly.

(2) Suppose that $M=\oplus_{i\in\frac 1 2\Z_+}M(i)$ is an
admissible $V$-module. Then (a) $\Omega_n(M)\supset \oplus_{i\leq
n}M(i)$; (b) Assume that $M$ is simple. Then
$\Omega_n(M)=\oplus_{i\leq n}M(i)$, and each $M(i)$ is a simple
$A_n(V)$-module for $i=0,\frac 1 2,..., n$.

(3) $M \mapsto M(0)$ gives a bijection between irreducible admissible $V$-modules and simple $A(V)$-modules

(4) The identity map induces an epimorphism from $A_{n}(V)$ to
$A_{m}(V)$ for any $n\geq m.$

(5) If $V$ is $g$-rational there are only finitely many irreducible admissible $g$-twisted $V$-modules up to isomorphism and each irreducible admissible $g$-twisted $V$-module is ordinary.
\end{theo}

We should point out that part (3) of the theorem was obtained in
\cite{KW}.

The next lemma will be used as a characterization of rationality
of $V$ in terms of semisimplicity of $A_n(V)$ for large enough
$n$.

\begin{lemm}\label{generalizedvsapce}
Suppose that $A(V)$ is finite dimensional, then any admissible
$V$-module is a direct sum of generalized eigenspaces for $L(0)$.
\end{lemm}
\noindent{\bf{Proof.}}\ \  Let $M=\oplus_{i\in\frac 1 2\Z_+}M(i)$ be
an admissible $V$-module with $M(0)\neq 0$.  Let $W$ be a maximal
subspace of $M$ which is a direct sum of generalized eigenspaces
with respect to $L(0)$. Then it is not hard to see that $W$ is a
submodule of $M$. Consider the $A(V)$-module $M(0)$. By our
assumption on finite dimension of $A(V)$ we see that there exists a
nonzero simple $A(V)$-submodule of $M(0)$, on which $L(0)$  acts as
a scalar by Schur lemma. This shows that $W\neq 0$.  We shall show
$W=M$. Suppose $M/W\neq 0$. Choose the minimal $n\in\frac 1 2\Z_+$
such that $M(n)/W(n)\neq 0$, where $W(n)=W\cap M(n)$. Then by the
similar argument as above, we see that $M(n)/W(n)$ contains a
nonzero simple $A(V)$-submodule, say $\mathscr W(n)/W(n)\neq 0$,
where $\mathscr W(n)$ is a subspace of $M(n)$. Since both $\mathscr
W(n)/W(n)$ and $W(n)$ are  a direct sum of generalized eigensapces
for $L(0)$, then so is  $\mathscr W(n)$. Thus $\mathscr W(n)\subset
W$ and $\mathscr W(n)=W(n)$,  a contradiction.  \QED

\vskip5pt Assume that $A(V)$ is finite dimensional. Let
$f(x)=(x-\lambda_1)(x-\lambda_2)\cdots(x-\lambda_r)\in\C[x]$ be the
monic polynomial   of least degree such that $f([w])=0$ in $A(V)$.
Then on any given simple $A(V)$-module $L(0)$ must act as a constant
$\lambda_i$ for some $i$. Note from Theorem \ref{theo3.2} that $V$ has exactly
$r$ inequivalent irreducible admissible modules $M^i=\sum_{n\in\frac 1 2 \Z_+}M^i_{\lambda_i+n}$ for
$i=1,\cdots,r.$ Then there exists $m_i>0$ such that $M^i_{\lambda_i+n}\ne 0$ for all $n\geq m_i.$
Let $N$ be a positive integer greater than
$| \lambda_i-\lambda_j\mid,$ $|\lambda_i|+1,$ and $m_i$ for $i,j=1,...,r$.

Note that the rationality is defined from the representation theory.
It is always believed that such property, which is analogous to the
semisimplicity  of Lie and associative algebras,  should have its own internal
characterization. The the following result can be regraded as an internal characterization of rationality.

\begin{theo}\label{propsemisimple}
$V$ is rational if and only if $A_n(V)$ is finite dimensional and
semisimple for some $n\geq N$.
\end{theo}
\noindent{\bf{Proof.}}\ \ The proof of Theorem 4.10 in \cite{DLM0}
shows that $V$ is rational then $A_n(V)$ is semisimple and finite
dimensional for all $n.$  Now we assume that $A_n(V)$ is
semisimple for some $n\geq N.$ By Theorem \ref{theo3.2} $A_m(V)$
is semisimple for all $m\leq n.$  Let $M=\oplus_{i\in\frac 1
2\Z_+}M(i)$ be an admissible $V$-module with $M(0)\neq 0$. By
Lemma \ref{generalizedvsapce}, we can write
$$M=\sum_{\lambda\in\{\lambda_1,...,\lambda_r\}}\bigoplus_{n\in\frac{1}
{2}\Z_+}M_{\lambda+n}$$ where $M_{\lambda+n}$ is the generalized
eigenspace for $L(0)$ with eigenvalue $\lambda+n$. Note that for
each $\lambda\in \{\lambda_1,...,\lambda_r\}$ the subspace
$\bigoplus_{n\in\frac{1}{2}\Z_+}M_{\lambda+n}$ is an admissible
submodule of $M.$ Without loss of generality, we may assume that
$M=\oplus_{n\in\frac 1 2\Z_+}M(n)$ for some
$\lambda\in\{\lambda_1,...,\lambda_r \}$ where
$M(n)=M_{\lambda+n}$.

We assert that the submodule $W$ generated by $\oplus_{n\leq N,\
n\in\frac 1 2\Z_+}M(n)$ is equal to the entire $M$. Otherwise,
$0\neq M/W=\oplus_{n>N,\ n\in\frac 1 2\Z_+}M(n)/W(n)$ where
$W(n)=W\cap M(n).$  Let $n_0\in\frac 1 2\Z_+$ be minimal such that
$M(n_0)/W(n_0)\neq 0$. Then $n_0>N$ and $M(n_0)/W(n_0)$ is an
$A(V)$-module by Theorem \ref{theo3.2}. Since $A(V)$ is
semisimple, there exists a nonzero simple $A(V)$-submodule of
$M(n_0)/W(n_0)$ on which $L(0)$ acts as the constant
$\lambda+n_0\in\{\lambda_1,...,\lambda_r\}$, which implies
$|\lambda-\lambda_j|=n_0$ for some $j$. But this is impossible by
our choice on $N.$ Thus we must have $W=M$.

We next show that if $X$ is a simple $A(V)$-submodule then $X$
generates an irreducible $V$-module $U.$ Denote
$J=\oplus_{n\in\frac12\Z_+}J(n)$ the maximal submodule of $U$ such
that $J(0)=0$, where $J(n)=J\cap U(n)$. Then the quotient $W=U/J$
is irreducible and $W(0)=X.$ Since $\oplus_{0\leq n\leq N}U(n)$ is
a semisimple $A_N(V)$-module we can regard each $W(n)$ a
$A_N(V)$-submodule of $U(n)$ for $n\leq N.$ From the choice of $N$
we know that  $W(N)\ne 0.$ Then the admissible $V$-submodule of
$U$ generated by $W(N)$ contains $W(0)=X.$ Thus $W(N)=U(N)$ and
therefore $J(N)=0$. By our choice of $N$ again we see that $J$
must be trivial. This implies that $U=W$ is irreducible.

It follows that the admissible $V$-submodule $\mathcal W$ of $M$
generated by $M(0)$ is completely reducible. Note that
$M(1)=\mathcal W(1)\oplus P$ where $P$ is a semisimple
$A(V)$-module. Again the admissible submodule of $M$ generated by
$P$ is completely reducible. Continuing in this  way completes the
proof. \QED

\begin{rema} Even in the case that $V$ is a vertex operator algebra, Theorem \ref{propsemisimple} strengthens the Theorem 4.11 of \cite{DLM0} where we require that $A_n(V)$ is semisimple for all $n.$
\end{rema}

\begin{rema}\label{remark1} There is a twisted analogue $A_{g,n}(V)$ (cf. \cite{DLM4}) of $A_n(V).$
One can similarly define the positive integer $N_g$. Then Theorem
\ref{propsemisimple} still  holds. That is,  $V$ is $g$-rational
if and only if $A_{g,n}(V)$ is finite dimensional and semisimple
for some $n\geq N_g$.
\end{rema}

We now use Theorem \ref{propsemisimple} to prove the following
result:
\begin{prop}\label{proprational}
Let $V=V_{\bar 0}\oplus V_{\bar 1}$ be a VOSA. If $V_{\bar 0}$ is
rational, then $V$ is rational.
\end{prop}
\noindent{\bf{Proof.}}\ Suppose that $V_{\bar 0}$ is rational.
Then by Theorem \ref{propsemisimple},  $A_n(V_{\bar 0})$ is finite
dimensional semisimple associative algebra if $n$ is sufficiently
large. This implies that $A_n(V)$ is semisimple as $A_n(V)$ is a
quotient of $A_n(V_{\bar 0})$. Applying Theorem
\ref{propsemisimple} again yields that $V$ is rational. \QED
\bigskip

We remark that we do not know how to prove the rationality of $V$ from the rationality
of $V_{\bar 0}$ without using $A_n(V).$ It is certainly a very interesting problem to find
a different approach without using $A_n(V).$ Although we can not show the converse of  Proposition
\ref{proprational} we strongly believe that rationalities of $V$
and $V_{\bar 0}$ are equivalent.

In the rest of this section we use the extension functor to consider the
rationality of a vertex operator superalgebra $V$. This approach has been studied in \cite{A} for vertex operator algebra. But our rationality result is different from that given in \cite{A}.

First let us describe the set Ext$_V^1(M^2, M^1)$ for any weak $V$-module
$M^1$ and $M^2$. We call a weak $V$-module $M$ an {\em extension of
$M^2$ by $M^1$} if there is a short exact sequence $0\rightarrow
M^1\rightarrow M\rightarrow M^2\rightarrow 0$. Two extensions $M$
and $N$ of $M^2$ by $M^1$ are said to be equivalent if there exists
a $V$-homomorphism $f: M\rightarrow N$ such that the following
diagram commutes:
\begin{eqnarray*}
\begin{CD}
0@>>>M^1 @>\phi>> M@>\varphi>> M^2@>>>0\quad \rm(exact)\\
@. @| @VV f V @| @.\\
0@>>>M^1 @>\phi^\prime>> N@>\varphi^\prime>> M^2@>>>0\quad
\rm(exact).
\end{CD}
\end{eqnarray*}
Define Ext$_V^1(M^2,M^1)$ to be the set of all equivalent classes
of $M^2$ by $M^1$. It is well known that Ext$_V^1(M^2,M^1)$
carries the structure of an abelian group such that the equivalent
class of $M^1\oplus M^2$ is the zero element.

Here is another  equivalent condition of
rationality.

\begin{prop}
Let $V$ be a vertex operator superalgebra. Then $V$ is rational if
and only if the following two conditions hold:
\begin{itemize}
\item[\rm(a)] Every admissible $V$-module contains a nontrivial
irreducible admissible submodule;
\item[\rm(b)] For any irreducible $V$-modules $M$ and $N$,
$\mathrm{Ext}_V^1(M,N)=0$.
\end{itemize}
\end{prop}
\noindent{\bf Proof.} It is clear that rationality implies both (a) and (b). Now we assume that (a) and (b) hold.  Let
$M=\oplus_{n\in\frac12\Z_+}M(n)$ be a nonzero admissible $V$-module.
Let $W$ be the sum of irreducible admission $V$-submodules of $M.$
Then $W=\oplus_{i\in I}W^i$ where each $W^i$ is an irreducible
admissible $V$-module. By condition (a), $W\neq 0$. We assert that
$W=M$. Otherwise consider the quotient module $M/W$.
It follows from the condition (a) again that there exists a weak $V$-submodule
$M^\prime$ such that $M^\prime\varsupsetneq W$ and $M^\prime/W$ is
an irreducible admissible $V$-module. Then by condition (b) and the
properties of Ext we have
$$\mathrm{Ext}_{V}^1(M^\prime/W, W)=\bigoplus_{i\in
I}\mathrm{Ext}_{V}(M^\prime/W, W^i)=0.$$ That is, $M^\prime
=M^\prime/W\oplus W$ as $V$-modules, contradicting the maximality
of $W$. So the assertion is true and $M$ is a direct sum of
irreducible admissible $V$-modules.\QED
\bigskip

We now turn our attention to the generators of vertex operator
superalgebras.

\begin{prop}\label{propa}  Let $V$ be vertex operator superalgebra. Then we have

(a) If $V$ is rational or $C_2$-cofinite, then $V$ is finitely generated.

(b) If $V$ is finitely generated, then $Aut(V)$ is an algebraic group.
\end{prop}

The proof of these results have been given in the case of vertex operator algebras in \cite{DZ},
\cite{KL} (also see \cite{GN})  and \cite{DG}. The same proof works here.

\vskip10pt \noindent{\bf 4. Vertex operator subalgebra generated
by
$V_{\frac12}$}\setcounter{section}{4}\setcounter{theo}{0}\setcounter{equation}{0}
\vskip5pt

In this section we study the vertex operator super subalgebra $U$
of $V$ generated by $V_{\frac12}$ and decompose $V$ as a tensor
product $U\otimes U^c$ where $U$ is holomorphic in the sense that
$U$ is the only irreducible module for itself and $U^c$ whose
weight $\frac12$ subspace is $0$ is the commutant of $U$ in $V$.
This decomposition reduces the study of vertex operator
superalgebras to the study of vertex operator superalgebras whose
weight $\frac12$ subspaces are $0.$

Let $V$ be a simple vertex operator
superalgebra of strong CFT type. Then there is a unique invariant, symmetric and nondegenerate bilinear form $(\cdot,\cdot)$ such that
\begin{equation}\label{normalized}
(\mathbf{1},\mathbf{1})=\sqrt{-1}
\end{equation}
(see \cite{L1} and \cite{X}). Then for $u, v\in
V_{\frac 1 2}$ one has
\begin{equation}\label{bili}
u_0v=(u, v)\mathbf{1}
\end{equation}  and
\begin{equation}\label{commutation}
[u(m), v(n)]_+=(u,v)\delta_{m+n+1,0}.
 \end{equation}  Note that the
restriction of $(\cdot,\cdot)$ to $V_{\frac 1 2}$ is still
nondegenerate.  Let
$\{a^1,a^2,\cdots,a^l\}$ be an orthonormal basis of $V_{\frac12}$
with respect to the form $(\cdot,\cdot)$ where $l=\dim V_{\frac12}.$

Let $U$ the vertex super subalgebra of $V$ generated by $V_{\frac
1 2}$. Then using the equation (\ref{commutation}) we see that
\begin{equation*}
U=\mathrm{Span}\{u^1_{-n_1}u^2_{-n_2}\cdots u^r_{-n_r}\mathbf 1\mid
u^i\in V_{\frac 1 2},\ n_1\geq n_2\geq\cdots\geq n_r>0\
\mathrm{and}\ r\in\Z_+\}.
\end{equation*}  In fact, $U$ carries the structure of a vertex operator
 superalgebra with conformal vector
 $$
\omega^\prime=\frac 12\sum_{i=1}^la^i_{-2}a^i_{-1}\mathbf 1.$$
Define operators $L^\prime(n)$ for $n\in\Z$ by:
\begin{equation*}
Y(\omega^\prime,z)=\sum_{n\in\Z}w^\prime_nz^{-n-1}=\sum_{n\in\Z}L^\prime(n)z^{-n-2}.
\end{equation*}
Then the weight $n$ subspace $U_n$ for  $L^\prime(0)$ is given as
follows
$$U_n=\big\langle u^1_{-n_1}u^2_{-n_2}\cdots u^r_{-n_r}\mathbf 1\big|
\begin{array}{c}u^i\in V_{\frac 1 2},
\ n_1\geq n_2\geq\cdots\geq n_r>0,\ r\in\Z_+, \\
 \mathrm{and}\ n_1+n_2\cdots+n_r=n+\frac r 2
\end{array}\big\rangle.
$$

It is well known (cf. \cite{KW}) that the vertex operator algebra $U$
generated by $V_{\frac 1 2}$ is holomorphic. So for any
admissible $V$-module $M$, we can decompose $M$ into irreducible
$U$-modules as follows
\begin{equation*} M=U\otimes \bar M,
\end{equation*}
where $\bar M=\{w\in M|u_nw=0\ \mathrm{for\ all}\ u\in U$ and
$n\in\Z_+\}$ is the multiplicity space of $U$ in $M.$ If $M=V$ the
multiplicity space $\bar M$ is denoted by $U^c$ and is called the
commutant of $U$ in $V$. In particular,  $V=U\otimes U^c.$ The
$U^c$ is a vertex operator superalgebra (see \cite{FZ} and
\cite{LL}) with $\omega-\omega^\prime$ as its conformal vector and
$U^c_{\frac 1 2}=0.$

 Let $Irr(V)$ and  $Irr(U^c)$ denote the sets of the isomorphism
classes of admssible irreducible $V$-modules and
$U^c$-modules, respectively. The following result is straightforward.
\begin{prop} Let $V$ be a vertex operator superalgebra. Then
\begin{itemize}\parskip-3pt
\item[\rm(a)\it]For any admissible  $V$-module $M,$ $\bar M$ is an admissible $U^c$-module. Moreover,  $M$ is irreducible if and
only if $\bar M$ is irreducible;
\item[\rm(b)\it] The map $U\otimes*:\ Irr(U^c)\rightarrow Irr(V)$ is a bijection.
\item[\rm(c)\it] $V$ is rational if and only if $U^c$ is rational.
\end{itemize}
\end{prop}

\vskip10pt\noindent{\bf 5.  The structure of weight 1 subspace}
\setcounter{section}{5}\setcounter{theo}{0}\setcounter{equation}{0}
\vskip5pt In this section we will investigate the Lie algebra
structure of weight 1 subspace $V_1$ and show that $V_1$ is a
reductive Lie algebra if $V$ is $\sigma$-rational using the
modular invariance results obtained in \cite{DZ1}. We also find an
upper bound for the rank of $V_1$ in terms of effective central
charge. Similar results for vertex operator algebra were given
previously in \cite{DM} and the proof presented here is a
modification of that used in \cite{DM}. We also apply these
results to estimate the dimension of weight $\frac12$ subspace
$V_{\frac12}$ of $V.$

We need some discussion on vertex operator superalgebra on torus (\cite{Z} and \cite{DZ2}), vector-valued modular forms \cite{KM} and the modular invariance of trace functions (\cite{Z} and \cite{DZ2}).

Let $V$ be a vertex operator superalgebra. The vertex operator superalgebra
$$(V, Y[v,z], {\bf 1}, \tilde\omega)$$
 on torus (see \cite{Z} and  \cite{DZ2}) is defined as follows:
\begin{eqnarray*}&Y[v,z]=Y(v, e^z-1)e^{z\mathrm{wt}
v}=\sum\limits_{n\in\Z}v[n]z^{-n-1}\\
&Y[\tilde w, z]=\sum_{n\in\Z}L[n]z^{-n-2} \end{eqnarray*} for
homogeneous $v$ and $\tilde \omega=\omega-\frac{c}{24}$.

We denote the eigenspace of $L[0]$ with eigenvalue $n\in\frac12\Z$
by $V_{[n]}$. If $v\in V_{[n]}$ we write wt$[v]=n$.

A holomorphic {\em vector-valued modular form of weight k} ($k$ any real number)
on the modular group $\Gamma=SL(2,\Z)$ may be described as follows:
for any integer $p\geq 1$ it is a tuple $(T_1(\tau),...,T_p(\tau))$
of functions holomorphic in the complex upper half-plane
$\mathbf{H}$ together with a $p$-dimensional complex representation
$\rho:\ \Gamma\rightarrow$ GL($p,\C$) satisfying the following
conditions:
\begin{itemize}\parskip-3pt
\item[(a)] For all $\gamma\in\Gamma$ we have
\begin{equation*}
(T_1,...,T_p)^t|_k
\gamma(\tau)=\rho(\gamma)(T_1(\tau),...,T_p(\tau))^t
\end{equation*}
(t refers to transpose of vectors and matrices).
\item[(b)] Each function $T_j(\tau)$ has a convergent $q$-expansion
holomorphic at infinity:
\begin{equation*}
T_j(\tau)=\sum_{n\geq 0}a_n(j)q^{n/N_j}
\end{equation*}
for positive integer $N_j$. (Here and below, $q$=exp$2\pi i\tau$).
\end{itemize}
The following result \cite{KM} plays an important role in this section.
\begin{prop}\label{vector-valued}
Let $(T_1,...,T_p)$ be a holomorphic vector-valued modular form of
weight $k$ associated to a representation $\rho$ of $\Gamma$. Then
there is a nonnegative constant $\alpha$ depending only on $\rho$
such that the Fourier coefficients $a_n(j)$ satisfy the polynomial
growth condition $a_n(j)=O(n^{k+2\alpha})$ for every $1\leq j\leq
p$.
\end{prop}

Fix automorphisms $g,h$ of $V$ of finite orders. Let $M$ be a
simple $g\sigma$-twisted $V$-module. Then $M=\oplus_{n=0}^\infty
M_{\lambda+\frac {n} {T^\prime}}$ for some $\lambda$ called the
conformal weight of $M$ ($M_\lambda\neq 0$), where $T^\prime$ is
the order of $g\sigma$. Suppose that $M$ is $\sigma h$-stable,
which is equivalent to the existence of a linearly isomorphic map
$\phi(\sigma h):M\rightarrow M$ such that
$$\phi(\sigma h)Y_M(v,z)\phi(\sigma h)^{-1}=Y_M((\sigma h)v,z)$$
for all $v\in V$. From now on we assume that $V$ is $C_2$-cofinite. Then function $F_M$ which is
linear in $v\in V$ is defined for homogenous $v\in V$ as follows:
\begin{equation}\label{tracef}
F_M(v,\tau)=q^{\lambda-c/24}\sum\limits_{n=0}^\infty\mathrm{tr}_{M_{\lambda+\frac
nT}}o(v)\phi(\sigma h)q^{\frac n
T}=\mathrm{tr}_Mo(v)\phi(\sigma h)q^{L(0)-c/24}
\end{equation}
which is a holomorphic function in the upper half plane
$\mathbf{H}$ \cite{DZ2}. Here and below we write $F_M(\tau)$
rather than $F_M(\mathbf{1},\tau)$ for simplicity. Then for any
$u,v$ in $V$ such that $gv=hv=v$ we have
\begin{equation}\label{identity}
\mathrm{tr}_{M}o(u)o(v)\phi(\sigma h)q^{L(0)-c/{24}}=F_M(u[-1]v,\tau)
-\sum_{k\geq1}E_{2k}(q)F_M(u[2k-1]v,\tau)
\end{equation}
(see \cite{DZ2} and \cite{Z}).
The functions $E_{2k}(\tau)$ are the Eisenstein series of
weight $2k:$
$$
E_{2k}(q)=\frac{-B_{2k}}{2k!}+\frac
2{(2k-1)!}\sum\limits_{n=1}^\infty\sigma_{2k-1}(n)q^n$$
where  $\sigma_k(n)$ is the sum of the $k$-powers of the divisors of
$n$ and $B_{2k}$ is a Bernoulli number. The
$E_2(\tau)$  enjoys an exceptional transformation law. Namely, its
transformation with respect to the  matrix $S=\left(
                                               \begin{array}{cc}
                                                 0 & -1 \\
                                                 1 & 0 \\
                                               \end{array}
                                             \right)
$ as follows:
\begin{equation}\label{stranformation}
E_2(\frac{-1}{\tau})=\tau^2E_2(\tau)-\frac{\tau}{2\pi i}.
\end{equation}

We also need results on $1$-point functions on torus from \cite{DZ2}. Let $g,h$ be
be automorphisms of $V$ of finite orders. The space of $(g,h)$ $\mathrm{1}$-point functions $\mathcal
{C}(g,h)$ is the $\C$-linear space consisting of functions
$$S:\ V\times \mathbf{H}\rightarrow \C$$
satisfying certain conditions (see \cite{DZ2} for details). The following results can be found in \cite{DZ2}.
\begin{theo}\label{theo6.2} Let $V$ be $C_2$-cofinite and $g,h\in Aut(V)$ of finite orders. Then
(1) For $S\in\mathcal{C}(g,h)$ and $\gamma=\left(
                                         \begin{array}{cc}
                                           a & b \\
                                           c & d \\
                                         \end{array}
                                       \right)\in\Gamma,$
we define
$$
S|_\gamma(v,\tau)=S|_k(v,\tau)=(c\tau+d)^{-k}S(v,\tau)$$
for $v\in V_{[k]}$, and extend linearly. Then
$S|_\gamma\in\mathcal{C}((g,h)\gamma))$, where $(g,h)\gamma=(g^ah^c,
g^bh^d)$.

(2)  Let $M$ be a simple $g\sigma$-twisted $V$-module such that
$M$ is $h$ and $\sigma$-stable. Then $F_M(v,\tau)\in\mathcal{C}(g,h)$.

(3) Suppose that $V$ is $g\sigma $-rational and $M^1,...,M^m$ are
the inequivalent, simple $g\sigma $-twisted $V$-module such that
$M^i$ is $h$ and $\sigma$-stable. Let $F_1,...,F_m$ be the
corresponding trace functions defined by
$\mathrm{(\ref{tracef})}$. Then $F_1,...,F_m$ form a basis of
$\mathcal{C}(g,h)$.

\end{theo}

We now assume that $V$ is of strong CFT type. Recall from
\cite{FHL} that the weight 1 subspace $V_1$ of $V$  carries a
natural Lie algebra structure, the Lie bracket being given by
$[u,v]=u_0v$ for $u, v\in V_1$. Then any weak $V$-module is
automatically a $V_1$-module such that $v\in V_1$ acts as $v_0.$
Note that there is a unique nondegenerate symmetric invariant
bilinear form $\langle\cdot,\cdot\rangle$ such that
$\langle\mathbf 1,\mathbf 1\rangle=-1$ and the restriction of
$\langle\cdot,\cdot\rangle$ to $V_1$ endows $V_1$ with a
nondegenerate, symmetric, invariant bilinear form such that
$u_1v=\langle u, v\rangle\mathbf 1$ for $u,v\in V_1$.

The following two theorems are extensions of similar results from vertex operator algebras \cite{DM} to vertex operator superalgebras.
\begin{theo}\label{reductive}
Let $V$ be a strongly rational or strongly $\sigma$-rational. Then the Lie algebra $V_{1}$ is reductive.
\end{theo}

\noindent{\bf {Proof.}}  We first deal with the case that $V$ is $\sigma$-rational. We have to show that the nilpotent
radical $N$ of the Lie algebra $V_1$ is zero. Suppose not, take any nonzero element $u\in N$.
Each $V_i$ for $i\in\frac12\Z$ is finite dimensional $V_1$-module and has a composition series
$0=W^0\subset W^1\subset W^2\subset W^3\subset\cdots\subset$ such that  $u_0$ acts trivially on each
composite factor $W^i/W^{i-1}(i=1,2,\cdots)$.
Note that we can take $\phi(\sigma)=\sigma$ on $V.$
Thus $V$ is $\sigma$-stable. In fact, any
irreducible $V$-module is $\sigma$-stable (see Lemma 6.1 of \cite{DZ2}).
As a result, $tr_{V_i}o(u)o(v)\sigma=0$ for all $v\in V_1$ and $i\in\frac12\Z.$ It follows from (\ref{identity}) that
\begin{equation}\label{vannish}
F_V(u[-1]v,\tau)=\sum_{k\geq1}E_{2k}(\tau)F_V(u[2k-1]v,\tau)
\end{equation}
where $(g,h)=(\sigma,1)$ and $F_V\in {\cal C}(\sigma,1)$ by
Theorem \ref{theo6.2}.

Note that if $k>1$ is an integer, then the element $u[2k-1]v$ has
$L[0]$-weight $2-2k<0$ and hence is $0.$  The non-degeneracy of
the bilinear form $\langle\cdot,\cdot\rangle$ guarantees that
there exists $v\in V_1$ such that $\langle u,v\rangle=1$. With
this choice of $v$, (\ref{vannish}) simplifies to read
\begin{equation}\label{simplify}
F_V(u[-1]v,\tau)=E_{2}(\tau)F_V(\tau).
\end{equation}

By Theorem \ref{theo3.2}, $V$ has finitely many irreducible $\sigma$-twisted $V$-modules up to isomorphism. We denotes these modules by $M^1,...,M^m.$ Note from Theorem \ref{theo6.2} that the $S=
\left(\begin{array}{cc}0 & 1 \\ -1 & 0 \end{array}\right)\in\Gamma$ maps ${\cal C}(\sigma,1)$ to
${\cal C}(1,\sigma).$ By Theorem \ref{theo6.2} again we see that
$$F_V(u[-1]v,-\frac1\tau)=\tau^2\sum_{i=1}^ms_iF_{M^i}(u[-1]v,\tau)$$
and
$$F_V(-\frac1\tau)=\sum_{i=1}^ms_iF_{M^i}(\tau)$$
for some $s_i\in \C.$
  Similar to
equality (\ref{simplify}) we also have
\begin{equation*}\label{simplifym}
F_{M^i}(u[-1]v,\tau)=E_2(\tau)F_{M^i}(\tau)\ \  \mathrm{for}\ i=1,...,m.
\end{equation*}
Thus
\begin{eqnarray*}
\tau^2\sum_{i=1}^ms_iF_{M^i}(u[-1]v,\tau)&=&F_V(u[-1]v,\frac{-1}{\tau})\nonumber\\
                                     &=&E_2(\frac{-1}{\tau})F_V(\frac{-1}{\tau})\nonumber\\
                                     &=&(\tau^2E_2(\tau)-\frac{\tau}{2\pi
                                     i})\sum_{i=1}^ms_iF_{M^i}(\tau)\\
                                     &=&\tau^2\sum_{i=1}^ms_iF_{M^i}(u[-1]v,\tau)-\frac{\tau}{2\pi
                                     i}\sum_{i=1}^ms_iF_{M^i}(\tau)\nonumber.
\end{eqnarray*}
Of course, the equality (\ref{stranformation}) is involved in the
calculations above. Canceling the term
$\tau^2\sum_{i=1}^ms_iF_{M^i}(u[-1]v,\tau)$ gives rise to the identity
$\sum_{i=1}^ms_iF_{M^i}(\tau)=0$, which in turn implies
$F_V(\frac{-1}\tau)=0$. But this is clearly not true, since
$$F_V(\frac{-1}\tau)=q^{-\frac{c}{24}}(\sum_{n\in \Z}(\dim V_n)q^n-\sum_{n\in \frac 12+\Z}(\dim V_n)q^n)\neq 0.$$
 So $N=0$ and  $V_1$ is
reductive.

Now we assume that $V$ is rational. As before we need to show that the nilpotent
radical $N$ of the Lie algebra $V_1$ is zero. This time we use
${\cal C}(\sigma,\sigma)$ instead of ${\cal C}(\sigma,1)$ and ${\cal C}(1,\sigma).$ In this case, the
$S\in \Gamma$ maps ${\cal C}(\sigma,\sigma)$ to itself. The similar argument just applies. \QED

\begin{rema} It is proved in \cite{DM} that if a vertex operator algebra is a strongly rational then weight one subspace is reductive. If one can prove that the rationality of $V_{\bar 0}$ from the rationality and $\sigma$-rationality of $V$ then Theorem \ref{reductive} will follows immediately. Unfortunately, none of these results have been established.
\end{rema}

The following result will be used in the next section.
\begin{lemm}\label{lemaa} Let $V$ be a vertex operator superalgebra.

(a) If $V$ is strongly rational, then any admissible
 $V$-module is completely reducible $V_1$-module. This is also equivalent to saying the action of any Cartan subalgebra of Lie algebra $V_1$ is semisimple  on any admissible $V$-modules.

 (b) If $V$ is strongly $\sigma$-rational, then any admissible
 $\sigma$-twisted $V$-module is completely reducible $V_1$-module.

 (c) If $V$ is either strongly rational or strongly $\sigma$-rational, then any irreducible $\sigma^i$-twisted $V$-module is an semisimple $V_1$-module for $i=0,1.$
\end{lemm}

\noindent{\bf{Proof.}} Since the proof of (b) is similar to (a) we
only show (a) and (c) for strongly rational vertex operator
superalgebra $V.$ Let $H$ be a Cartan subalgebra of $V_1.$ It is
enough to show that $H$ acts semisimply on any irreducible
$\sigma^i$-twisted $V$-module for $i=0,1.$ Since the homogeneous
subspaces of an irreducible $\sigma^i$-twisted $V$-module is
always finite dimensional there is a common eigenvector of $H$ on
the irreducible module. So it is enough to show that $H$ acts on
$V$ semisimply.

First we show that for any $0\ne u\in H,$ $h_0$ is not nilpotent. Note that the restriction of the bilinear form $\langle\cdot,\cdot\rangle$ to $H$ is nondegenerate. If $u_0$ for some $0\ne u\in H$ is nilpotent, we can take $v\in H$ such that $\langle u,v\rangle=1.$  The proof of Theorem \ref{reductive} then gives a contradiction.

We now prove that $u_0$ is semisimple on $V.$ Since $Aut(V)$ is an algebraic group by Proposition \ref{propa} and $\{e^{tu_0}|t\in\C\}$ is one dimensional algebraic subgroup of $Aut(V)$ we see immediately
see that $\{e^{tu_0}|t\in\C\}$ is isomorphic to the 1-dimensional multiplicative algebraic group $\C_m$
as $u_0$ is not nilpotent (cf. \cite{Ma}). This finishes the proof. \QED
\bigskip

Now that $V_1$ is reductive, there are two extreme cases: $V_1$ is
a semisimple Lie algebra and $V_1$ is abelian. The vertex operator
subalgebra generated by $V_1$ will be extensively investigated in
Section 6. We next study the rank of $V_1$ in the rest of this
section. Let $l$ be the rank of $V_1.$ That is, $l$ is the
dimension of a Cartan subalgebra $H$ of $V_1.$  Similar to the
case of  vertex operator algebras  in \cite{DM}, $l$ is closed
related to the {\em effective central charge} $\tilde c$ which is
defined as follows: Let $\{M^1,...,M^m\}$ be the irreducible
$\sigma$-twisted $V$-modules up to isomorphism. Then there exist
$\lambda_i\in \C$ such that
$M^i=\sum_{n\in\frac12\Z_+}M^i_{\lambda_i+n}$ with
$M^i_{\lambda_i}\ne 0.$ The $\lambda_i$ is called the conformal
weight of $M^i.$ By Theorem 8.9 of \cite{DZ2}, $\lambda_i$ and the
central charge $c$ of $V$ are rational numbers for all $i.$ Define
$\lambda_{\mathrm min}$ to be the minimum of the conformal weights
$\lambda_i$ and set
\begin{equation*}\label{effectivecharge}
\tilde c=c-24\lambda_{\mathrm min},\ \ \
\tilde\lambda_i=\lambda_i-\lambda_{\mathrm min}.
\end{equation*}

\begin{theo}\label{dimcartan}
Let $V$ be strongly $\sigma$-rational. Then $l\leq \tilde c$
\end{theo}

\noindent{\bf{Proof.}} Let $H$ be a Cartan  subalgebra of
$V_1$. Note that the component operators of the vertex operators
$Y(u,z)$ on $V$ for $u\in H$ form a Heisenberg Lie algebra. This
amounts to saying that for $u,v\in H$ the following relations hold:
\begin{equation}
[u_m,v_n]=m\delta_{m,-n}\langle u,v\rangle.
\end{equation}
In fact, these relations also hold true on any $\sigma$-twisted
$V$-module $M.$

Consider $(g,h)=(1,1).$ Let $F_i=F_{M^i}$ be as defined in
(\ref{tracef}). Then
 $F_i\in {\cal C}(1,1).$  Recall
 that
 $$\eta(\tau)=q^{1/24}\prod_{n=1}^\infty(1-q^n)$$
 is a modular form of weight $\frac12.$  Then
\begin{equation*}\eta(\tau)^{\tilde c}F_i(\tau)=q^{\tilde
\lambda_i}\prod_{n=1}^\infty(1-q^n)^{\tilde
c}\sum\limits_{n=0}^\infty\mathrm{tr}_{M_{\lambda_i+\frac
n{\mathrm2}}}\phi(\sigma )q^{\frac n {\mathrm2}}
\end{equation*}
is holomorphic in $\mathbf H\cup \{i\infty\}$.  Now it follows from
the transformation law for $\eta(\tau)$ and Theorem \ref{theo6.2}
that the $m$-tuple
\begin{equation*}
\big(\eta(\tau)^{\tilde c}F_1(\tau),...,\eta(\tau)^{\tilde
c}F_m(\tau)\big)
\end{equation*}
is a holomorphic vector-valued modular form of weight ${\tilde c}/2$. So the
Fourier coefficients of $\eta(\tau)^{\tilde c}F_i(\tau)$ have
polynomial growth by Proposition \ref{vector-valued}.

 The Stone-von-Neumann theorem provides us a somewhat different way
 to look more closely at $F_i(\tau)$. Namely,
$M^i$ has the following tensor decomposition:
\begin{equation}\label{tensordec}
M^i=M(1)\otimes_{\C}\Omega_{M^i},
\end{equation}
where $M(1)=\C[u_m|u\in H,m>0]$ is the Heisenberg vertex operator
algebra of rank $l$ generated by $H$ and $\Omega_{M^i}=\{w\in
M^i|u_nw=0$ for $u\in H$ and $n>0$\}. Then the trace function
$F_i(\tau)$ corresponding to the decomposition (\ref{tensordec}) is
equal to
\begin{equation*}
q^{(l-c)/24}\eta(\tau)^{-l}\mathrm{tr}_{\Omega_i}\phi(\sigma)q^{L(0)},
\end{equation*}
as tr$_{M(1)}\phi(\sigma)q^{L(0)}=q^{l/24}\eta(\tau)^{-l}$. Thus
\begin{equation}\label{alternative}
\eta(\tau)^{\tilde c}F_i(\tau)=q^{(l-c)/24}\eta(\tau)^{\tilde
c-l}\mathrm{tr}_{\Omega_i}\phi(\sigma)q^{L(0)}.
\end{equation}
We know that the  Fourier coefficients of the left-hand side
of (\ref{alternative}) have polynomial growth. This forces the same
is true on $\eta(\tau)^{\tilde c-l}$. Then one has
$\tilde c-l\geq 0$, as $\eta(\tau)^s$ has exponential growth of
Fourier coefficients whenever $s<0$ (cf. \cite{K}).\QED

\vskip3pt We now use Theorem
\ref{dimcartan} to do an estimation on the dimension of
$V_{\frac 12}$.

\begin{coro}\label{main} Let $V$ be strongly $\sigma$-rational. Then
$\mathrm{dim}V_{\frac 1 2}\leq 2\tilde c+1$.
\end{coro}

\noindent{\bf{Proof.}}  Let $d$ be a nonnegative integer such that
$2d\leq$dim$V_{\frac 1 2}\leq 2d+1$. Then there exists a unique
(up to a constant) nondegenerate bilinear form satisfying
(\ref{normalized}).  Note that the restriction of $(\cdot,\cdot)$
to $V_{\frac 1 2}$ is still nondegenerate.  So we can choose
elements $b^i, {b^i}^*\in V_{\frac 1 2}$ such that $(b^i,
{b^j}^*)=\delta_{ij}$ and $(b^i, b^j)=0=({b^i}^*,{b^j}^*)$ for all
$1\leq i,j\leq d$.  Set $h^i=b^i_{-1}(b^i)^*_{-1}{\bf 1}$ for
$i=1,...,d.$ Then $h^i\in V_1$ and $h^i_1h^j=\delta_{i,j},$
$h^i_0h^j=0$ for $i,j\in\{1,\cdots,d\}.$ As a result,
$\sum_{i=1}^d\C h^i\subset V_1$  is contained in a Cartan
subalgebra of $V_1.$ By Theorem \ref{main}, $d\leq l\leq \tilde c$
and the proof is complete. \QED

\vskip10pt \noindent{\bf 6. $C_2$-confiniteness and
Integrability}\setcounter{section}{6}\setcounter{theo}{0}\setcounter{equation}{0}
\vskip5pt

We continue our discussion on  the weight 1 subspace $V_1.$ We
will determine the vertex operator subalgebra $\langle V_1\rangle$
of $V$ generated by $V_1$ following the approach in \cite{DM1}. It
turns out that $\langle V_1\rangle$ is isomorphic to $L_{\mathfrak
g_1}(k_1,0)\otimes \cdots \otimes L_{\mathfrak g_s}(k_s,0)\otimes
M(1)$ where $V_1={\mathfrak g_1}\oplus \cdots \oplus{\mathfrak
g_s}\oplus Z(V_1),$ ${\mathfrak g_i}$ are simple, $k_i\geq 1$ are
integers and $M(1)$ is the Heisenberg vertex operator algebra
built up from $Z(V_1)$ (see below for the definition of
$L_{\mathfrak g}(k,0)$). Moreover, $\langle V_1\rangle$ is
contained in the rational vertex operator subalgebra $L_{\mathfrak
g_1}(k_1,0)\otimes \cdots \otimes L_{\mathfrak g_s}(k_s,0)\otimes
V_L$ for some positive definite lattice $L\subset Z(V_1)$
satisfying $rank(L)=\dim Z(V_1).$

Here we need to review the construction of untwisted affine
Kac-Moody Lie algebras $\widehat{\mathfrak g}$ associated with
simple Lie algebras $\mathfrak g$ and relevant results from \cite{Kac}. Let $\mathfrak h$
be a Cartan subalgebra  of $\mathfrak g$ and $\Phi$ the corresponding root system.
Fix a nondegenerate symmetric invariant bilinear form $(\cdot,\cdot)$ on $\widehat{\mathfrak g}$ such that
the square length of a long root is $2$ where we have identified $\mathfrak h$ with its dual via the bilinear form. Then the affine Kac-Moody algebra associated to $\mathfrak g$ is given by
$$
\widehat{\mathfrak g}=\mathfrak g\otimes\C[t,t^{-1}]\oplus\C K$$
with the bracket relations
\begin{equation}\label{lieb}
[u(m),v(n)]=[u,v](m+n)+m(u,v)\delta_{m+n,0}K,\  [K,\widehat{\mathfrak g}]=0
\end{equation}
for $u,v\in\mathfrak g$ and $m,n\in\Z$ where $u(m)=u\otimes t^m$.
Let $L(\lambda)$ be the irreducible $\mathfrak g$-module with
highest weight $\lambda\in\mathfrak h$. Consider $L(\lambda)$ as a
$\mathfrak g\otimes\C[t]$-module with $\mathfrak g\otimes\C t[t]$
acting trivially and with $K$ acting as the scalar  $k\in\C$. Then
the generalized Verma module
$$V(k,\lambda)=\mathrm{Ind}_{\mathfrak g}^{\widehat{\mathfrak
g}}L(\lambda)=U(\widehat{\mathfrak g})\otimes_{U(\mathfrak
g\otimes\C [t]\oplus\C K)}L(\lambda)$$ has the unique irreducible
quotient $L(k,\lambda).$ It is well known  that $L(k,\lambda)$ is
integrable if, and only if, $k$ is a nonnegative integer and
$\lambda$ is a dominant integral weight such that
$(\lambda,\theta)\leq k$ where $\theta\in \Phi$ is the maximal
root.

Let $V$ be a VOSA of strong CFT type and $\langle\cdot,\cdot\rangle$
be the unique nondegenerate bilinear form satisfying $\langle\mathbf
1,\mathbf 1\rangle=-1.$ Suppose that $\mathfrak g\subset V_1$ is a simple subalgebra. Then both  bilinear forms $(\cdot,\cdot)$ and
$\langle\cdot,\cdot\rangle$ on $\mathfrak g$ are symmetric and
invariant, so they must be proportional, that is,
\begin{equation}\label{defk}
\langle \cdot,\cdot\rangle=k(\cdot,\cdot)\ \ \ \mathrm{for\ some}\
k\in\C.
\end{equation}
Then for any $u,v\in V_1$ and integers $m,n$ one has
\begin{equation*}
[u_m,v_n]=[u,v]_{m+n}+mu_1v\delta_{m+n,0}.
\end{equation*}
Comparing this with (\ref{lieb}) shows that  the map
$$u(m)\rightarrow u_m\ \ \mathrm{for}\ u\in\mathfrak g\ \mathrm{and}\
m\in\Z$$
together with $K\rightarrow k$ gives rise to a representation of
$\widehat{\mathfrak g}$ of level $k$.

Now we are going to state our  main result related to
$C_2$-integrability, which has already been proved to be true in
\cite{DM1} for vertex operator algebras satisfying the
$C_2$-confiniteness. But given a vertex operator superalgebra
$V=V_{\bar 0}\oplus V_{\bar 1}$ which satisfies the $C_2$-cofinite
condition, generally we can not prove that the even part  $V_{\bar
0}$ also has such property. So in this sense, the following result
sharpens the Theorem 3.1 of \cite{DM1} although the idea is similar.
\begin{theo}\label{simplepart}
Let $V$ be a simple vertex operator superalgebra which is
$C_2$-cofinite of strong CFT type, with $\mathfrak g\subset V_1 $
a simple Lie subalgebra, $k$ the level of $V$ as
$\widehat{\mathfrak{g}}$-module, and  the vertex operator
subalgebra $U$ of $V$ generated by $\mathfrak g$. Then the
following hold:\parskip-2pt
\begin{itemize}\parskip-3pt
\item[$\mathrm(a)$] The restriction of $\langle\cdot,\cdot\rangle$ to $\mathfrak g$ is
nondegenerate,
\item[$\mathrm{(b)}$] $U\cong L(k,0)$,
\item[$\mathrm{(c)}$] $k$ is a positive integer,
\item[$\mathrm(d)$] $V$ is an integrable $\widehat{\mathfrak g}$-module.
\end{itemize}
\end{theo}

\noindent{\bf{Proof.}} Let $\mathfrak h$ be a Cartan subalgebra of
$\mathfrak g$ and let $\mathfrak g=\mathfrak
h\oplus\sum_{\alpha\in\Phi}\mathfrak g_\alpha$ be the corresponding
Cartan decomposition of $\mathfrak g$. Since $\mathfrak g$ is generated by subalgebras isomorphic to
$sl(2,\C)$ it is good enough to show the theorem for  $\mathfrak g=sl(2,\C).$ Let
$\{h,x,y\}$ be the standard basis of $\mathfrak g.$  Then $(\alpha,\alpha)=2$ and
  $k=\frac{\langle\alpha,\alpha\rangle}2$ from this and equation (\ref{defk}).

Clearly, $U=\langle\mathfrak g\rangle$ is a quotient of $V(k,0)$. So $U$ is a
$\widehat{\mathfrak g}$-integrable module if and only if $\mathcal
U=L(k,0)$ for some $k\in\Z_+.$  This is also equivalent to the
existence of a positive integer $r$ such that
\begin{equation}\label{killed}
(x_{-1})^r\mathbf1=0.
\end{equation}
The proof of (\ref{killed}) is similar to the same result in \cite{DM1} and we omit the proof.
(b) then follows immediately. Also note that
$\mathfrak g\subset U$, so $ U$ can not be a
one-dimensional  trivial module. Thus $k\neq 0$ and $k$ must be a
positive integer, proving (c) and (a). Since
$L(k,0)$ is rational (cf \cite{DLM1}),  $V$ is a direct sum of
irreducible $L(k,0)$-modules, each of which is integrable as
$\widehat{\mathfrak g}$-module. Hence $V$ is an integrable
$\widehat{\mathfrak g}$-module. This proves (d). \QED
\vskip3pt

Next we consider a toral subalgebra of $V_1.$  Let $V$ be strongly rational or strongly $\sigma$-rational  and $\mathfrak h\subset V_1$ be a toral subalgebra such that the
restriction of $\langle\cdot,\cdot\rangle$ to $\mathfrak h$ remains
nondegenerate. Notably, any Cartan subalgebra of $V_1$ automatically
satisfies such condition.

\begin{theo}\label{abelianpart}
Suppose that $V$ is strongly rational or strongly $\sigma$-rational. Let
$\mathfrak h\subset V_1$ be a toral subalgebra such that the
restriction of $\langle\cdot,\cdot\rangle$ to $\mathfrak h$ is
nondegenerate. Then there exist a positive-definite even lattice
$L\subset\mathfrak h$ with rank {\rm dim}$\mathfrak h$  and a vertex
operator supersubalgebra $U$ of $V$ such that $\mathfrak
h\subset U\cong V_L$.
\end{theo}

This Theorem has been proved in \cite{DM1} (also see \cite{Ma}) for vertex operator algebra. The same argument using Lemma \ref{lemaa} is also valid for vertex operator superalgebra.

We now assume that
$$V_1={\mathfrak g_1}\oplus\cdots \oplus {\mathfrak g_s}\oplus {\mathfrak h}$$
where ${\mathfrak g_i}$ are simple Lie algebras and $Z(V_1)={\mathfrak h}.$ By Theorems \ref{reductive}, \ref{simplepart} and \ref{abelianpart} we have (see \cite{DM1} and \cite{Ma}):
\begin{coro} The  $V$ contains a strongly rational vertex operator subalgebra
$$U=L_{\mathfrak g_1}(k_1,0)\otimes\cdots \otimes L_{\mathfrak g_s}(k_s,0)\otimes V_L$$
where and the commutant $U^c$ of $U$ in $V$ is a vertex operator superalgebra such that
$U^c_1=0.$
 \end{coro}


\vskip5pt


\begin{thebibliography}{AAAA}

\bibitem{A}T. Abe, Rationality of the vertex operator algebra $V_L^+$ for a positive definite even lattice $L,$  {\em Math. Z.} {\bf  249} (2005),  455–484.
\bibitem{ABD}T. Abe, G. Buhl and  C. Dong,  Rationality, regularity, and $C_2$-cofiniteness,
 {\em Trans. Amer. Math. Soc.} {\bf 356} (2004), 3391-3402.

\bibitem {B} R. E. Borcherds, Vertex algebras, Kac-Moody algebras,
and the Monster, {\em Proc. Natl. Acad. Sci. USA} {\bf 83} (1986), 3068-3071.

\bibitem{DVVV} R. Dijkgraaf, C. Vafa, E. Verlinde and H.
Verlinde, The operator algebra of orbifold models, {\em Comm. Math.
Phys.} {\bf 123} (1989), 485-526.


\bibitem{DG}C. Dong and R. Griess Jr., Automorphism groups and derivation algebras of finitely generated vertex operator algebras, {\em Michigan Math. J.}  {\bf 50} (2002), 227–239.

\bibitem{DL}C. Dong and J. Lepowsky, Generalized vertex algebras and
relative vertex operators, Progress in Mathematics, 112.
Birkh\"{a}user Boston, Inc., Boston, MA, 1993.

\bibitem{DLM1} C. Dong, H. Li and G. Mason,
Regularity of rational vertex operator algebras, {\em  Adv. Math.} {\bf 132} (1997), 148-166.

\bibitem{DLM2} C. Dong, H. Li and  G. Mason,  Twisted representations of
vertex operator algebras, {\em  Math. Ann.} {\bf 310} (1998), 571-600.

\bibitem{DLM0}C. Dong, H. Li and G. Mason,  Vertex operator algebras and associative
algebras, {\em J Algebra} {\bf  206}  (1998), 67-96.

\bibitem{DLM4}C. Dong,  H. Li and  G. Mason,  Twisted representations
of vertex operator algebras and associative algebras, {\em Internat.
Math. Res. Notices} {\bf 8} (1998), 389-397.

\bibitem{DLM3} C. Dong, H. Li and G. Mason, Modular invariance of trace
functions in orbifold theory and generalized moonshine, {\em Comm.
Math. Phys.} {\bf 214} (2000), 1-56.

\bibitem{DM}C. Dong and G. Mason, Rational vertex operator
algebras and the effective central charge, {\em Int. Math. Res. Not.} {\bf 56} (2004),
2989-3008.

\bibitem{DM1} C. Dong and  G. Mason, Integrability of
$C_2$-cofinite vertex operator algebras, {\em Int. Math. Res. Not.} (2006), Art. ID 80468.

\bibitem{DZ} C. Dong and W. Zhang, Rational vertex operator algebras are finitely generated, {\em  J. Algebra} {\bf 320} (2008),  2610–2614.

\bibitem{DZ2} C. Dong and  Z. Zhao, Modularity in orbiford theory for vertex operator
superalgebras, {\em Comm. Math. Phys.} {\bf 260} (2005), 227-256.

\bibitem{DZ1}C. Dong and Z. Zhao,  Twisted representations of vertex operator superalgebras,
 {\em Commun. Contemp. Math.} {\bf 8} (2006), 101-121.

\bibitem{FFR} A. J. Feingold, I. B. Frenkel and  John F. X. Ries, Spinor Construction
of Vertex Operator Algebras, Triality, and $E_{8}^{(1)}$, {\em Contemp. Math.} {\bf 121} (1991).

\bibitem{FHL} I. B. Frenkel, Y. Huang and J. Lepowsky,  On axiomatic approaches to
 vertex operator algebras and modules, {\em Mem. Amer. Math. Soc.} {\bf 104} (1993).


\bibitem{FLM2} I. B. Frenkel, J. Lepowsky and A. Meurman,
Vertex Operator Algebras and the Monster, {\em Pure and Applied Math.}
Vol. {\bf 134}, Academic Press, 1988.

\bibitem{FZ}I. B. Frenkel and Y. Zhu, Vertex operator algebras associated to
representations of affine and Virasoro algebras, {\em Duke Math. J.} {\bf 66}
(1992),  123-168.

\bibitem{GN} M. Gaberdiel and A. Neitzke, Rationality,
quasirationality and finite $W$-algebras,  {\em Comm. Math. Phys.} {\bf 238}
(2003),  305-331.

\bibitem{Kac}V. G.  Kac, Infinite-dimensional Lie algebras, Third
edition, Cambridge University Press, Cambridge, 1990.

\bibitem{KW} V. Kac and W. Wang, Vertex operator superalgebras and their
representations, {\em  Contemp. Math.} {\bf 175} (1994), 161-191.

\bibitem{KL} M. Karel and H. Li,  Certain generating subspaces for vertex operator algebras,
 {\em J. Algebra} {\bf 217} (1999), 393-421.

\bibitem{K} M. Knopp, Modular Functions in Analytic Number Theory, Markham
Publishing, Illinois, 1970.

\bibitem{KM}M. Knopp and G. Mason,
On vector-valued modular forms and their Fourier coefficients,  {\em Acta
Arith.} {\bf 110} (2003), 117--124.


\bibitem{LL}J. Lepowsky and H. Li, Introduction to vertex operator algebras
and their representations, {\em Progress in Mathematics} {\bf 227}
Birkh\"{a}user Boston, Inc., Boston, MA, 2004.

\bibitem{L1} H. Li, Symmetric invariant bilinear forms on vertex operator
 algebras, {\em J. Pure Appl. Algebra} {\bf 96} (1994), 279-297.

\bibitem{L2} H. Li, Local systems of vertex operators, vertex superalgebras and
modules, {\em J. Pure Appl. Algebra} {\bf 109} (1996), 143-195.

\bibitem{L3} H. Li, Local systems of twisted vertex operators; vertex superalgebras and twisted
modules, {\em Contemp. Math.} {\bf 193} (1996), 203-236.

\bibitem{L4} H. Li, On abelian coset generalized vertex algebras, {\em Commun. Contemp. Math.} {\bf 3} (2001), 287-340.

\bibitem{Ma} G. Mason, Lattice subalgebras of strongly regular vertex operator
algebras, arXiv:1110.0544.


\bibitem{X}X. Xu, Introduction to Vertex Operator Superalgebras and Their Modules, {\em  Mathematics and its Applications} Vol. {\bf 456},
 Kluwer Academic Publishers, Dordrecht, 1998.

\bibitem{Z} Y. Zhu, Modular invariance of characters of vertex operator algebras, {\em J. Amer, Math. Soc.} { \bf 9} (1996), 237-302.

\end{thebibliography}
\end{document}